\documentclass[12pt]{amsart}
\usepackage{amsthm, amstext, amsmath, amssymb, graphicx, verbatim}
\usepackage[
  left=1in,
  right=1in
]{geometry}
\usepackage{booktabs}
\usepackage{hyperref}
\usepackage{ mathrsfs }
\usepackage{float}
\usepackage{algorithm}
\usepackage[dvipsnames]{xcolor}
\usepackage{versions}
\usepackage{enumerate, enumitem}
\usepackage{graphicx}
\usepackage{mathtools}

\usepackage{algorithm}
\usepackage{tabularx}
\usepackage{MnSymbol}
\usepackage{subcaption}
\usepackage[normalem]{ulem} 
\usepackage{tikz-cd}
\usepackage{colonequals}

\usepackage[utf8]{inputenc}
\usepackage[english]{babel}
\DeclareFixedFont{\ttb}{T1}{txtt}{bx}{n}{12} 
\DeclareFixedFont{\ttm}{T1}{txtt}{m}{n}{12}

\newtheorem{thm}{Theorem}[section]
 
\newtheorem{lemma}[thm]{Lemma} 
\newtheorem{prop}[thm]{Proposition}
\theoremstyle{definition}
\newtheorem{defn}[thm]{Definition}

\newtheorem{remark}[thm]{Remark}

\newcommand{\R}{\mathbb R}
\newcommand{\C}{\mathbb C}

\newcommand{\abs}[1]{\left|#1\right|}

\setenumerate[0]{label=(\alph*)}

 \allowdisplaybreaks

\DeclareMathOperator{\Disc}{Disc}

\begin{document}
\raggedbottom
 
\title{The dynamics of cyclic snap on polygons}


\author[Alvarez]{Nicolas Jose Alvarez} 
\address{Department of Mathematics, Elon University, Elon, NC, 27244}
\email{nalvarez2@elon.edu}

\author[Chen]{Jiaxuan Chen}
\address{Stanford Online High School, Redwood City, CA 94063}
\email{jiaxchen@ohs.stanford.edu}

\author[Dietrich]{Anna Dietrich}
\address{Department of Mathematics, Brown University, Providence, RI, 02912}
\email{anna\_dietrich@brown.edu}

\author[Miller]{Steven J. Miller}
\address{Department of Mathematics, Williams College, Williamstown , MA 01267}
\email{sjm1@williams.edu}

\author[Lee]{Uihyeon Lee}
\address{Department of Mathematics, Texas A\&M University, College Station, TX 77840}
\email{uihyeonlee@tamu.edu}

\author[Weinreich]{Max Weinreich}
\address{Department of Mathematics, CUNY Baruch College, New York, NY 10010}
\email{max.weinreich@baruch.cuny.edu}

\author[Zhao]{Chenzixi Zhao}
\address{Department of Mathematics, Purdue University, West Lafayette, IN 47907}
\email{zhao1452@purdue.edu}

\begingroup

\renewcommand{\thefootnote}{*}

\footnotetext{Research of this paper is partially supported by the U.S. National Science Foundation under Grant \#DMS2341670.}

\endgroup

\begin{abstract}
 We examine the limiting behavior of polygon cyclic snapping, an operation inspired by the midpoint iteration. We show that every $N$-gon converges under cyclic snap to a weighted sum of the initial vertices and that the rate of convergence slows asymptotically as the number of vertices grows.
\end{abstract}

\maketitle

\section{Introduction}

We introduce \emph{cyclic snap}, a linear dynamical system acting on polygons inspired by the classical midpoint iteration of Darboux \cite{darboux, MR3541977, MR3986100}. We locate the collapse point of polygon snap and study the rate of convergence using linear algebra and polynomial techniques.

\begin{defn}
A $N$-gon is an $N$-tuple $v = (v_0, \ldots, v_{N-1})$ of points in $\R^2$. We interpret indices modulo $n$, so $v_n = v_0$. We view $v$ as an element of the vector space $(\R^2)^N$ and call $v_i$ the $i$-th vertex of $v$. For each $0 \leq i \leq N - 1$, we define the \emph{snap} operation at vertex $i$ by
$$s_i: (\mathbb{R}^2)^N\to (\mathbb{R}^2)^N,$$
$$s_i(v_0,\dots,v_i,\dots,v_{N-1}) \ = \ (v_0,\dots,v_i', \dots v_{N-1}),$$
where $v'_i$ is the midpoint of $v_{i-1}$ and $v_{i+1}$. The \emph{cyclic snap} operation is defined by
$$T : (\mathbb{R}^2)^N\to (\mathbb{R}^2)^N, $$
$$T \ =\  s_{N-1} \circ \ldots \circ s_0.$$
\end{defn}
Intuitively, the operation $s_i$ flattens the $i$-th corner of the polygon by moving that vertex to the midpoint of its neighboring vertices. Whereas the midpoint iteration applies all these flattening operations simultaneously, cyclic snap applies them in order.
    
\begin{thm} \label{thm_main}
    Let $N \geq 1$ and let $v = (v_0, \ldots, v_{N-1})$ be an $N$-gon. Then the vertices of the iterates $v, T(v), T^2(v),\dots$ of $v$ under the cyclic snap operation converge to the point in $\R^2$ given by the weighted average $\frac{2v_{N-1}}{N} + \frac{v_{N-2}}{N} + \dots + \frac{v_1}{N}$.
\end{thm}

\begin{figure}[tbp]
\centering
\includegraphics[width=\linewidth]{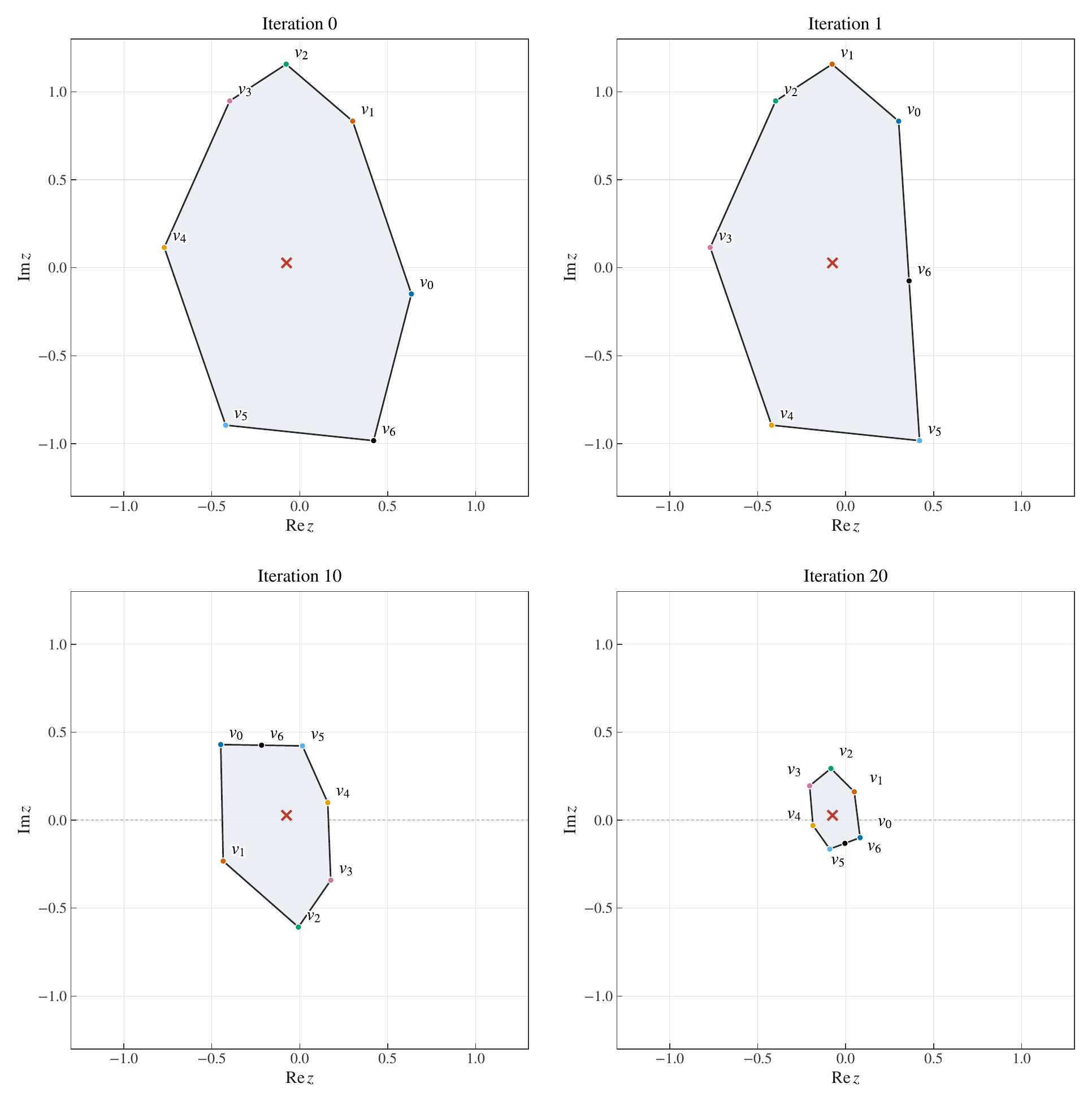}
\caption{A convex heptagon after $0$, $1$, $10$, and $20$ successive snap updates, generated with our applet \cite{polymath_polygon_dynamics}. Each update moves one vertex to the midpoint of its two neighbors, proceeding cyclically; seven updates constitute one application of $T$. All panels use the same coordinate scale. The red cross marks the common limit point given by Theorem~\ref{thm_main}, and the rate of collapse to the limit point is described by Proposition \ref{prop_rate}.}
\label{fig:heptagon-snap}
\end{figure}

Figure~\ref{fig:heptagon-snap} illustrates the contraction of a convex heptagon toward the limit point in Theorem~\ref{thm_main}.

Our study belongs to the general area of polygon dynamics, a field that combines aspects of algebra, geometry and dynamical systems. The objects of study in polygon dynamics are self-mappings of spaces of polygons. A prototypical example is the midpoint iteration of Darboux, which sends an input polygon (given as a list of vertices) to the polygon formed by the midpoints of the edges in cyclic order \cite{darboux}. Iterating this operation produces a sequence of polygons that converges to the centroid of the initial input. In other words, the midpoint iteration asymptotically collapses a polygon to a limit point, namely, the centroid. Cyclic snap is a variant of the midpoint iteration in which the snap operations at each vertex are applied simultaneously rather than sequentially. We observe that, unlike the midpoint iteration, cyclic snap forgets the data of the $0$-th vertex (when $n > 1$) and therefore does not converge to the centroid in general.

\begin{remark}
 We were motivated by the study of polygon dynamics in $\R^2$ and therefore state our results only in that context. However, results generalize in straightforward fashion to polygons in $\R^m$ for any $m$; indeed, our proof proceeds by studying the coordinatewise evolution of the vertices, so it is perhaps most naturally stated in $\R^1$. Similarly, our main results extend to polygons in $\C^m \cong \R^{2m}$. In contrast, the cyclic snap operation does not necessarily converge over general valued fields. For example, the polygons diverge $2$-adically, since increasing powers of $2$ appear as denominators.
\end{remark}

Other polygon iterations, such as the pentagram map, also feature collapsing behavior. For nonlinear dynamical systems, proving existence of collapse points and finding exact formulas for them is an active field of study \cite{MR4403964, glick, stiegler}. Other polygon iterations may not have collapse points; for example, quadrilateral refolding is a system that can produce drift or quasiperiodic motion on an elliptic curve \cite{MR2135187,esch2001screensaver}.

In Section 2, we illustrate the proof for quadrilaterals; in this case, all the computations are totally explicit. In Section 3, we prove Theorem \ref{thm_main}. In Section 4, we propose some future research directions. In Section 5, we briefly discuss the code used in this project.

This research project took place as a part of the Polymath Jr. 2026 Program. The main results were proved in July-August 2026. The authors thank the organizers of the program for providing the opportunity for our collaboration.

\subsection*{AI Disclosure}
AI was used only for standard tasks, and was not used to create any conjectures or proofs. In particular, it was used by some of the students for the following purposes only: debugging code, improving code formatting, implementing the ``randomize" function for $n$-gons, and literature search. Anything generated by AI was checked before used.

\section{Quadrilaterals}

\begin{prop} \label{prop_quad}
Let $v = (v_0, v_1, v_2, v_3)$ be a $4$-gon. Then $v$ converges under the cyclic snap operation to $\frac{v_3}{2} + \frac{v_2+v_1}{4}$.
\end{prop}
\begin{proof}
    Let $x, y$ be the coordinates on $\R^2$. We prove convergence coordinatewise. By symmetry, it suffices to just check the claim for the $x$-coordinate, i.e., for all $0 \leq i \leq 3$, we have
        $$\lim_{k \to \infty} x(T^k(v)_i) \ =\  \frac{x(v_3)}{2} + \frac{x(v_2)+x(v_1)}{4}.$$
    
    Let $v$ be a quadrilateral with vertices $v_0, v_1, v_2, $ and $v_3$. Define a sequence $(a_k)_{k \geq -2}$ in $\R^2$ by 
    $$a_{-2} = x(v_1), \quad a_{-1} = x(v_2), \quad a_{0} = x(v_3), \quad a_1 = x(v_0'), \quad a_2 = x(v_1'), \quad a_3 = x(v_2') \dots. $$
    We claim that 
    \begin{equation}
        \lim_{k \to \infty} a_k \ = \  \frac{x(v_3)}{2} + \frac{x(v_2)+x(v_1)}{4}. \label{eq_ak_quad}
    \end{equation}
    This proves the proposition, since for each $0 \leq i \leq 3$, the sequence $(a_k)_{k \geq -2}$ contains $(x(T^j(v)_i))_{j \geq 0}$ as a subsequence.
    
    The sequence $a_k$ is a linear recurrence given for all $k \geq 1$ by $a_k = \frac{a_{k-1} + a_{k-3}}{2}$. As a matrix equation,
    $$\begin{bmatrix}
        a_k \\
        a_{k-1} \\
        a_{k-2}
    \end{bmatrix}
    \ =\ 
    \begin{bmatrix}
        1/2 & 0 & 1/2 \\
        1 & 0 & 0 \\
        0 & 1  & 0 
    \end{bmatrix}
    \begin{bmatrix}
        a_{k-1} \\
        a_{k-2} \\
        a_{k-3}
    \end{bmatrix}
    \ = \  
    \begin{bmatrix}
        1/2 & 0 & 1/2 \\
        1 & 0 & 0 \\
        0 & 1  & 0 
    \end{bmatrix}^2
    \begin{bmatrix}
        a_{k-2} \\
        a_{k-3} \\
        a_{k-4}
    \end{bmatrix}
    \ = \  \dots \ = \  
    \begin{bmatrix}
        1/2 & 0 & 1/2 \\
        1 & 0 & 0 \\
        0 & 1  & 0 
    \end{bmatrix}^{k}
    \begin{bmatrix}
        a_{0} \\
        a_{-1} \\
        a_{-2}
    \end{bmatrix}.$$
    Let 
    $$ A \ = \  
    \begin{bmatrix}
        1/2 & 0 & 1/2 \\
        1 & 0 & 0 \\
        0 & 1 & 0  
    \end{bmatrix}.$$
    Diagonalizing $A$ allows us to compute the matrix powers
    $$A^k \ = \  
    \begin{bmatrix}
    1 & \frac{-3+i\sqrt{7}}{8} & \frac{-3 - i\sqrt{7}}{8} \\
    1 & \frac{-1-i\sqrt{7}}{4} & \frac{-1 + i\sqrt{7}}{4} \\
    1 & 1 & 1
    \end{bmatrix}
    \cdot
    \begin{bmatrix}
    1 & 0 & 0 \\
    0 & \frac{-1-i\sqrt{7}}{4} & 0 \\
    0 & 0 & \frac{-1+i\sqrt{7}}{4}
    \end{bmatrix}^k
    \cdot
    \begin{bmatrix}
    \frac{1}{2} & \frac{1}{4} & \frac{1}{4} \\
    \frac{-7-5i\sqrt{7}}{28} & \frac{-7+ 11i\sqrt{7}}{56} & \frac{21-i\sqrt{7}}{56} \\
    \frac{-7+5i\sqrt{7}}{28} & \frac{-7-11i\sqrt{7}}{56} & \frac{21+i\sqrt{7}}{56}
    \end{bmatrix} $$
    where 
    $\lambda_1 = 1, \lambda_2 = \frac{-1 - i\sqrt{7}}{4}, \lambda_3 = \frac{-1 + i\sqrt{7}}{4} $ are the eigenvalues of $A$.
    Observing that $\abs{\lambda_2}, \abs{\lambda_3} < 1$, we obtain
    {\small
    \begin{align*}
    \lim\limits_{k \to \infty}
    \begin{bmatrix}
        a_k \\
        a_{k-1} \\
        a_{k-2} 
    \end{bmatrix} 
    \ &= \  
    \lim_{k \to \infty} A^k
    \begin{bmatrix}
        x(v_3) \\
        x(v_2) \\
        x(v_1) 
    \end{bmatrix} \\ 
    &= \  
    \lim_{k \to \infty}\begin{bmatrix}
    1 & \frac{-3+i\sqrt{7}}{8} & \frac{-3 - i\sqrt{7}}{8} \\
    1 & \frac{-1-i\sqrt{7}}{4} & \frac{-1 + i\sqrt{7}}{4} \\
    1 & 1 & 1
    \end{bmatrix}
    \cdot
    \begin{bmatrix}
    1 & 0 & 0 \\
    0 & \frac{-1-i\sqrt{7}}{4} & 0 \\
    0 & 0 & \frac{-1+i\sqrt{7}}{4}
    \end{bmatrix}^k
    \cdot
    \begin{bmatrix}
    \frac{1}{2} & \frac{1}{4} & \frac{1}{4} \\
    \frac{-7-5i\sqrt{7}}{28} & \frac{-7+ 11i\sqrt{7}}{56} & \frac{21-i\sqrt{7}}{56} \\
    \frac{-7+5i\sqrt{7}}{28} & \frac{-7-11i\sqrt{7}}{56} & \frac{21+i\sqrt{7}}{56}
    \end{bmatrix} \cdot \begin{bmatrix}
        x(v_3) \\
        x(v_2) \\
        x(v_1) 
    \end{bmatrix} \\ 
    &= \  \begin{bmatrix}
    1 & \frac{-3+i\sqrt{7}}{8} & \frac{-3 - i\sqrt{7}}{8} \\
    1 & \frac{-1-i\sqrt{7}}{4} & \frac{-1 + i\sqrt{7}}{4} \\
    1 & 1 & 1
    \end{bmatrix}
    \cdot
    \begin{bmatrix}
    1 & 0 & 0 \\
    0 & 0 & 0 \\
    0 & 0 & 0
    \end{bmatrix}
    \cdot
    \begin{bmatrix}
    \frac{1}{2} & \frac{1}{4} & \frac{1}{4} \\
    \frac{-7-5i\sqrt{7}}{28} & \frac{-7+ 11i\sqrt{7}}{56} & \frac{21-i\sqrt{7}}{56} \\
    \frac{-7+5i\sqrt{7}}{28} & \frac{-7-11i\sqrt{7}}{56} & \frac{21+i\sqrt{7}}{56}
    \end{bmatrix} \cdot \begin{bmatrix}
        x(v_3) \\
        x(v_2) \\
        x(v_1) 
    \end{bmatrix} 
    \\ &=
    \ 
    \begin{bmatrix}
    1 & 0 & 0 \\
    1 & 0 & 0 \\
    1 & 0 & 0
    \end{bmatrix} \cdot \begin{bmatrix}
    \frac{1}{2} & \frac{1}{4} & \frac{1}{4} \\
    \frac{-7-5i\sqrt{7}}{28} & \frac{-7+ 11i\sqrt{7}}{56} & \frac{21-i\sqrt{7}}{56} \\
    \frac{-7+5i\sqrt{7}}{28} & \frac{-7-11i\sqrt{7}}{56} & \frac{21+i\sqrt{7}}{56}
    \end{bmatrix} \cdot \begin{bmatrix}
        x(v_3) \\
        x(v_2) \\
        x(v_1) 
    \end{bmatrix} 
    \\ 
    &= \ 
    \begin{bmatrix}
    \frac{1}{2} & \frac{1}{4} & \frac{1}{4} \\
    \frac{1}{2} & \frac{1}{4} & \frac{1}{4} \\
    \frac{1}{2} & \frac{1}{4} & \frac{1}{4}
    \end{bmatrix}
    \begin{bmatrix}
        x(v_3) \\
        x(v_2) \\
        x(v_1) 
    \end{bmatrix} 
    \\
    &= \ 
    \begin{bmatrix}
        \frac{x(v_3)}{2}+ \frac{x(v_2)+x(v_1)}{4} \\
        \frac{x(v_3)}{2}+ \frac{x(v_2)+x(v_1)}{4} \\
        \frac{x(v_3)}{2}+ \frac{x(v_2)+x(v_1)}{4} \\
    \end{bmatrix}.
    \end{align*}
    }
    The first row is the desired equation \eqref{eq_ak_quad}.
\end{proof}

\begin{remark}
        If a polygon $w$ is output by a snap operation $s_i$, then $w$ is supported on a (possibly degenerate) triangle, with the $i$-th vertex located at the midpoint of a side of the triangle. Applying $s_{i+1}$ bisects the supporting triangle, so upon each further application of snap, the area halves.
\end{remark}

\section{Convergence of any $N$-gon under the cyclic snap operation}

In this section, we prove Theorem \ref{thm_main} and study the rate of convergence. We will require the following elementary theorem from the theory of polynomials \cite{enestrom1920remarque}. The form we state is the most standard version; see \cite[Theorem 1.3]{MR3329717}.

\begin{thm}[Eneström–Kakeya Theorem]\label{kakeya}
    Let $P(z) = \sum_{k=0}^{n} a_k z^k = a_nz^n + a_{n-1}z^{n-1} + \dots + a_1 z + a_0$ be a polynomial of degree $n$ with real coefficients. If the coefficients satisfy
        $$0 \le a_0 \le a_1 \le \dots \le a_{n-1} \le a_n,$$
    then all complex zeros of $P(z)$ satisfy $\abs{z} \le 1.$
\end{thm}

 The key idea of the proof of Theorem \ref{thm_main}, just as in the proof for quadrilaterals, is to study a matrix $A$ that describes vertex evolution under cyclic snap. In the quadrilateral case, direct inspection showed that $A$ has a simple eigenvalue $1$ that dominates all the other eigenvalues, which allowed us to take the limit of $A^k$ as $k \to \infty$. The proof of this for general $N$ relies on the following technical lemma.

 \begin{lemma}\label{strict_inequality} 
    Let $N \geq 3$. If $\lambda$ is a complex zero of $P(\lambda) = 2\lambda^{N-2} + \lambda^{N-3} + \dots + \lambda^2 + \lambda + 1$, then $\lvert \lambda \rvert < 1.$ 
    \begin{proof}
    Since the coefficients of $\lambda$ are decreasing nonnegative real numbers, the complex zeros of $P(\lambda)$ each satisfy $\abs{\lambda} \leq 1$ by Theorem \ref{kakeya}. Assume by way of contradiction that there is a zero $\lambda \in \C$ of $P(\lambda)$ such that $\abs{\lambda} = 1$.
    Define a new polynomial 
        $$Q(\lambda) \ = \  (\lambda - 1)P(\lambda) \ = \  2\lambda^{N-1} - \lambda^{N-2} - 1.$$
    Then $Q(\lambda) = 0$. Now, it suffices to address two cases. 
    
    First, consider the case $\lambda^{N-2} = 1$. Then $\lambda^{N-2} \cdot \lambda = \lambda^{N-1} = \lambda.$ So $0 = Q(\lambda)  = 2\lambda -2$, so $\lambda = 1$. But $P(1) = N \neq 0$, a contradiction. 
    
    In the remaining case $\lambda^{N-2} \neq 1$, we may 
        rearrange $Q(\lambda ) = 0$ to obtain
            $\abs{2\lambda^{N-1}} = \abs{\lambda^{N-2} + 1}.$
        Since we assumed $\abs{\lambda} = 1,$ it follows that $\abs{2\lambda^{N-1}} = 2 \abs{\lambda}^{N-1} = 2.$ However, $\abs{\lambda^{N-2} + 1} < 2$ because $\lambda^{N-2} \neq 1.$ Therefore $\abs{2\lambda^{N-2}} \neq \abs{\lambda^{N-2} + 1}$, a contradiction.
        
    Thus, we have $\lvert \lambda \rvert < 1$ for any zero $\lambda$ as desired.    
\end{proof}
\end{lemma}

The proof of Theorem \ref{thm_main} also requires knowing that a certain polynomial has distinct roots, which we check using the polynomial discriminant.

\begin{lemma} \label{lem_distinct_roots}
    For each $N \geq 3$, the polynomial $Q(\lambda) = 2 \lambda^{N -1} - \lambda^{N-2} - 1$ has distinct roots.
\end{lemma}

\begin{proof}
    Let $Q^*(\lambda) \  = \ -\lambda^{N-1} Q(\lambda^{-1}) \ = \  \lambda^{N-1} + \lambda - 2$. Then $Q^*$ is the negative of the reciprocal polynomial of $Q$. Since $Q$ has constant term $-1$, it does not have $0$ as a root, and factoring each polynomial into linear factors over $\C$ shows that the roots of $Q^*$ are the reciprocals of the roots of $Q$. Therefore it suffices to show that the roots of $Q^*$ are distinct, i.e., that the discriminant $\Disc(Q^*)$ is nonzero. We use the well-known formula (see e.g. \cite[p. 621]{MR2286236} or \cite{discriminanttrinomial}) for the discriminant of a trinomial:
    $$\Disc( z^n + bz + c) \  = \ (-1)^{n(n-1)/2}  n^n c^{n-1} + (-1)^{(n-1)(n-2)/2} (n-1)^{n-1} b^n.$$
    Setting $n = N - 1$, $b = 1$, $c = -2$
    yields
    \begin{align*}
        \Disc( Q^*(\lambda) ) \  &= (-1)^{(N-1)(N-2)/2}  (N-1)^{N-1} (-2)^{N-2} + (-1)^{(N-2)(N-3)/2} (N-2)^{N-2} \\
        &= (-1)^{(N+1)(N-2)/2}  (N-1)^{N-1} 2^{N-2} + (-1)^{(N-2)(N-3)/2} (N-2)^{N-2} \\
        &= (-1)^{(N+1)(N-2)/2} \, \left[(N-1)^{N-1} 2^{N-2} + (N-2)^{N-2} \right].
    \end{align*}
    Since $N \geq 3$, the expression $(N-1)^{N-1} 2 ^{N-2} + (N - 2)^{N-2}$ is positive, so the discriminant is nonzero.
\end{proof}

\begin{proof}[Proof of Theorem \ref{thm_main}]
Let $v = (v_0, \ldots, v_{N-1})$ be an $N$-gon. If $N = 1$ or $N = 2$, the claim is immediate, since the sequences of vertices are constant after one step. So assume $N \geq 3$. As in the proof of Proposition \ref{prop_quad}, it suffices to prove coordinatewise convergence, and in particular convergence of the $x$-coordinates.  Define the sequence $(a_k)_{k \geq -N+1}$ by $x(T^j(v)_i) = a_{-N+1+jN+i}$ for all $j \geq 0, i \geq 0$. The first few terms are
$$a_{-N+1} = x(v_0), \quad a_{-N+2} = x(v_1), \quad \ldots, \quad a_{0} = x(v_{N-1}), \quad a_1 = x(T(v)_0), \ldots. $$
Then for all $k \ge 1$, we obtain the following recursive formula
    $$a_k \ = \  \frac{a_{k-1}+a_{k-N+1}}{2}$$
and its matrix representation 
\begin{align*}
\begin{bmatrix}
    a_k \\
    a_{k-1} \\
    a_{k-2} \\
    \vdots \\
    a_{k-N+2}
\end{bmatrix} \  &= \  
\begin{bmatrix}
    1/2 & 0 & 0 & \cdots & 1/2 \\
    1 & 0 & 0 & \cdots & 0 \\
    0 & 1 & 0 & \cdots & 0 \\
    \vdots & \vdots & \vdots & \ddots & \vdots \\
    0 & 0 & 0 & \cdots & 0
\end{bmatrix} \cdot
\begin{bmatrix}
    a_{k-1} \\
    a_{k-2} \\
    a_{k-3} \\
    \vdots \\ 
    a_{k-N+1}
\end{bmatrix}  \\ \  &= \  
\begin{bmatrix}
    1/2 & 0 & 0 & \cdots & 1/2 \\
    1 & 0 & 0 & \cdots & 0 \\
    0 & 1 & 0 & \cdots & 0 \\
    \vdots & \vdots & \vdots & \ddots & \vdots \\
    0 & 0 & 0 & \cdots & 0
\end{bmatrix}^2 \cdot 
\begin{bmatrix}
    a_{k-2} \\
    a_{k-3} \\
    a_{k-4} \\
    \vdots \\
    a_{k-N}
\end{bmatrix} \\
\ &= \ \dots \ = \  
\begin{bmatrix}
    1/2 & 0 & 0 & \cdots & 1/2 \\
    1 & 0 & 0 & \cdots & 0 \\
    0 & 1 & 0 & \cdots & 0 \\
    \vdots & \vdots & \vdots & \ddots & \vdots \\
    0 & 0 & 0 & \cdots & 0
\end{bmatrix}^k \cdot
\begin{bmatrix}
    a_0 \\
    a_{-1} \\
    a_{-2} \\
    \vdots \\
    a_{-N+2}
\end{bmatrix}.
\end{align*}
The matrix
    $$A \ = \ 
        \begin{bmatrix}
            1/2 & 0 & 0 & \cdots & 1/2 \\
            1 & 0 & 0 & \cdots & 0 \\
            0 & 1 & 0 & \cdots & 0 \\
            \vdots & \vdots & \vdots & \ddots & \vdots \\
            0 & 0 & 0 & \cdots & 0
        \end{bmatrix}$$
is a companion matrix, that is, each nonzero entry is either in the first row or is a $1$ on the subdiagonal. We compute the characteristic polynomial by cofactor expansion of $\lambda I_n - A$ along the first row:    
    \begin{align*}
    \det (\lambda I_n - A) \ &= \ \lambda^{N-1} -       \frac12 \lambda^{N-2} - \frac12 \\
    & \ = \ \frac{1}{2}Q(\lambda) \\
    &\ = \  \frac{1}{2}(\lambda-1)\underbrace{(2\lambda^{N-2} + \lambda^{N-3} + \dots + \lambda + 1)}_{P(\lambda)}.
    \end{align*}
So the eigenvalues of $A$ are $\lambda_1 = 1$ and the roots $\lambda_2, \lambda_3, \dots, \lambda_{N-1}$ of $P(\lambda)$. Lemma \ref{strict_inequality} shows that that $\lvert \lambda_2 \rvert, \lvert \lambda_3 \rvert, \dots, \lvert \lambda_{N-1} \rvert  < 1$. Lemma \ref{lem_distinct_roots} shows that the roots $\lambda_i$ are all distinct. Since the complex eigenvalues of $A$ are distinct, the matrix $A$ is diagonalized by a matrix of eigenvectors, i.e.,
$$A = PDP^{-1}$$
where $D$ is the diagonal matrix with entries $\lambda_1, \ldots, \lambda_n$ and the $i$-th column of $P$ is a $\lambda_i$-eigenvector. If $\lambda$ is an eigenvalue, then the characteristic equation $\lambda^{N-1} - \frac{1}{2} \lambda^{N-2} - \frac{1}{2} = 0$ holds, so a direct check shows that $w_\lambda := (\lambda^{N-2} \; \ldots \; \lambda^0)^{\mathrm{T}}$ is an eigenvector:
\begin{align*}
    Aw_\lambda \  &= \  
        \begin{bmatrix}
            1/2 & 0 & \cdots & 1/2 \\
            1 & 0 & \cdots & 0 \\
            \vdots & \vdots & \ddots & \vdots \\
            0  & 0 & \cdots & 0
        \end{bmatrix} \begin{bmatrix}
            \lambda^{N-2} \\
            \lambda^{N-1} \\
            \vdots \\
            \lambda^0 \\
        \end{bmatrix} \\
        \ &= \  \begin{bmatrix}
            \frac{1}{2}\lambda^{N-2} + \frac{1}{2} \\
            \lambda^{N-2} \\
            \vdots \\
            \lambda^1 \\
        \end{bmatrix} \\
    \ &= \  \begin{bmatrix}
            \lambda^{N-1} \\
            \lambda^{N-2} \\
            \vdots \\
            \lambda^1 \\
        \end{bmatrix} \\
        \ &= \  \lambda w_\lambda.
\end{align*}
Therefore we may choose the diagonalizing matrix of eigenvectors to be the Vandermonde matrix
$$
P \ = \ 
\begin{bmatrix}
    \lambda_1^{N-2} & \lambda_2^{N-2} & \lambda_3^{N-2} & \dots & \lambda_{N-1}^{N-2} \\
    \lambda_1^{N-3} & \lambda_2^{N-3} & \lambda_3^{N-3} & \dots & \lambda_{N-1}^{N-3} \\
    \lambda_1^{N-4} & \lambda_2^{N-4} & \lambda_3^{N-4} & \dots & \lambda_{N-1}^{N-4} \\
    \vdots & \vdots & \vdots & \ddots & \vdots \\
    \lambda_1^{0} & \lambda_2^{0} & \lambda_3^{0} & \dots & \lambda_{N-1}^{0} \\
\end{bmatrix} \ = \ 
\begin{bmatrix}
    1 & \lambda_2^{N-2} & \lambda_3^{N-2} & \dots & \lambda_{N-1}^{N-2} \\
    1 & \lambda_2^{N-3} & \lambda_3^{N-3} & \dots & \lambda_{N-1}^{N-3} \\
    1 & \lambda_2^{N-4} & \lambda_3^{N-4} & \dots & \lambda_{N-1}^{N-4} \\
    \vdots & \vdots & \vdots & \ddots & \vdots \\
    1 & 1 & 1 & \dots & 1
\end{bmatrix}.
$$
It turns out that, to find $\lim_{k \to \infty} a_k,$ it suffices to find only the first row of $P^{-1}$. This is because of the specific form of the matrices in the equation $A^k = P D^k P^{-1}$, as we now explain. Since $A = PDP^{-1},$ multiplying by $P^{-1}$ gives us $P^{-1} A = DP^{-1}.$ Now let $r_1 = 
\begin{bmatrix}
    x_1 & x_2 & \dots & x_{N-1}
\end{bmatrix}$ be the first row of $P^{-1}.$ As we are only concerned with the first row, we obtain
$$r_1 A \ = \ \lambda_1 r_1 \ = \  1 \cdot r_1 \ = \  r_1,$$
$$
\begin{bmatrix}
    x_1 & x_2 & x_3 & \dots & x_{N-1}
\end{bmatrix} 
\begin{bmatrix}
    1/2 & 0 & 0 & \cdots & 1/2 \\
    1 & 0 & 0 & \cdots & 0 \\
    0 & 1 & 0 & \cdots & 0 \\
    \vdots & \vdots & \vdots & \ddots & \vdots \\
    0 & 0 & 0 & \cdots & 0
\end{bmatrix} \ = \  
\begin{bmatrix}
    x_1 & x_2 & x_3 & \dots & x_{N-1}
\end{bmatrix},
$$
$$
\begin{bmatrix}
    \frac12 x_1 + x_2 & x_3 & x_4 & \dots & \frac12 x_1
\end{bmatrix} \ = \  
\begin{bmatrix}
    x_1 & x_2 & x_3 & \dots & x_{N-1}
\end{bmatrix}.
$$
This gives us $\frac12 x_1 = x_2 = x_3 = \dots = x_{N-1}.$ Furthermore, $P^{-1}P = I$ gives us that $x_1 + x_2 + x_3 + \dots + x_{N-1} = 1$. Therefore we obtain $x_1 = \frac{2}{N}, \  x_2 = \frac{1}{N},\  x_3 = \frac{1}{N}, \  \dots, \  x_{N-1} = \frac{1}{N}$. Combining everything and using $*$ to denote irrelevant entries, we get:
\begin{align*}
    A \ &= \  PDP^{-1} \\
    & \ = \  
    \begin{bmatrix}
    1 & \lambda_2^{N-2} & \lambda_3^{N-2} & \dots & \lambda_{N-1}^{N-2} \\
    1 & \lambda_2^{N-3} & \lambda_3^{N-3} & \dots & \lambda_{N-1}^{N-3} \\
    1 & \lambda_2^{N-4} & \lambda_3^{N-4} & \dots & \lambda_{N-1}^{N-4} \\
    \vdots & \vdots & \vdots & \ddots & \vdots \\
    1 & 1 & 1 & \dots & 1
\end{bmatrix} \cdot 
\begin{bmatrix}
    1 & 0 &0& \dots & 0 \\
    0 & \lambda_2 & 0 & \dots & 0 \\
    0 & 0 & \lambda_3 & \dots & 0 \\
    \vdots & \vdots & \vdots & \ddots & \vdots \\
    0 &0 & 0 & \dots & \lambda_{N-1}
\end{bmatrix}^k \cdot 
\begin{bmatrix}
    \frac{2}{N} & \frac{1}{N} & \frac{1}{N} & \dots & \frac{1}{N} \\
    * & * & * & \ddots & * \\
        * & * & * & \ddots & * \\
    \vdots & \vdots & \vdots & \ddots & \vdots \\
        * & * & * & \ddots & * \\
\end{bmatrix}.
\end{align*}
Taking the limit,
\begin{align*}
    \lim\limits_{k \to \infty} A^k 
    & \ = \  
    \begin{bmatrix}
    1 & \lambda_2^{N-2} & \lambda_3^{N-2} & \dots & \lambda_{N-1}^{N-2} \\
    1 & \lambda_2^{N-3} & \lambda_3^{N-3} & \dots & \lambda_{N-1}^{N-3} \\
    1 & \lambda_2^{N-4} & \lambda_3^{N-4} & \dots & \lambda_{N-1}^{N-4} \\
    \vdots & \vdots & \vdots & \ddots & \vdots \\
    1 & 1 & 1 & \dots & 1
\end{bmatrix} \cdot 
\begin{bmatrix}
    1 & 0 &0& \dots & 0 \\
    0 & 0 & 0 & \dots & 0 \\
    0 & 0 & 0 & \dots & 0 \\
    \vdots & \vdots & \vdots & \ddots & \vdots \\
    0 &0 & 0 & \dots & 0
\end{bmatrix} \cdot 
\begin{bmatrix}
    \frac{2}{N} & \frac{1}{N} & \frac{1}{N} & \dots & \frac{1}{N} \\
    * & * & * & \ddots & * \\
        * & * & * & \ddots & * \\
    \vdots & \vdots & \vdots & \ddots & \vdots \\
        * & * & * & \ddots & * \\
\end{bmatrix} \\
\ &= \  
\begin{bmatrix}
    1 & 0 & 0 & \dots & 0 \\
    1 & 0 & 0 & \dots & 0 \\
    1 & 0 & 0 & \dots & 0 \\
    \vdots & \vdots & \vdots & \ddots & \vdots \\
    1 & 0 & 0 & \dots & 0 \\
\end{bmatrix} \cdot 
\begin{bmatrix}
    \frac{2}{N} & \frac{1}{N} & \frac{1}{N} & \dots & \frac{1}{N} \\
    * & * & * & \ddots & * \\
        * & * & * & \ddots & * \\
    \vdots & \vdots & \vdots & \ddots & \vdots \\
        * & * & * & \ddots & * \\
\end{bmatrix} \\
& \ = \  
\begin{bmatrix}
    \frac{2}{N} & \frac{1}{N} & \frac{1}{N} & \dots & \frac{1}{N} \\
    \frac{2}{N} & \frac{1}{N} & \frac{1}{N} & \dots & \frac{1}{N} \\
    \frac{2}{N} & \frac{1}{N} & \frac{1}{N} & \dots & \frac{1}{N} \\
    \vdots & \vdots & \vdots & \ddots & \vdots \\
    \frac{2}{N} & \frac{1}{N} & \frac{1}{N} & \dots & \frac{1}{N}
\end{bmatrix}. \end{align*}

Therefore
\begin{align*} 
\lim\limits_{k \to \infty}
\begin{bmatrix}
    a_k \\
    a_{k-1} \\
    a_{k-2} \\ 
    \vdots \\
    a_{k-N+2}
\end{bmatrix} \ &= \  
\lim\limits_{k \to \infty} A^k \cdot 
\begin{bmatrix}
    a_0 \\
    a_{-1} \\
    a_{-2} \\
    \vdots \\
    a_{-N+2}
\end{bmatrix} \ = \ 
\begin{bmatrix}
    \frac{2}{N} & \frac{1}{N} & \frac{1}{N} & \dots & \frac{1}{N} \\
    \frac{2}{N} & \frac{1}{N} & \frac{1}{N} & \dots & \frac{1}{N} \\
    \frac{2}{N} & \frac{1}{N} & \frac{1}{N} & \dots & \frac{1}{N} \\
    \vdots & \vdots & \vdots & \ddots & \vdots \\
    \frac{2}{N} & \frac{1}{N} & \frac{1}{N} & \dots & \frac{1}{N}
\end{bmatrix} \cdot \begin{bmatrix}
    a_0 \\
    a_{-1} \\
    a_{-2} \\
    \vdots \\
    a_{-N+2}
\end{bmatrix} \\ \ &= \  
\begin{bmatrix}
    \frac{2a_0}{N} + \sum_{i = -N +2}^{-1} \frac{a_i}{N} \\
    \frac{2a_0}{N} + \sum_{i = -N +2}^{-1} \frac{a_i}{N} \\
     \frac{2a_0}{N} + \sum_{i = -N +2}^{-1} \frac{a_i}{N} \\
    \vdots \\
     \frac{2a_0}{N} + \sum_{i = -N +2}^{-1} \frac{a_i}{N}
\end{bmatrix} \ = \ 
\begin{bmatrix} 
    x(\frac{2v_{N-1}}{N} + \sum_{i=1}^{N-2} \frac{v_i}{N}) \\ 
    x(\frac{2v_{N-1}}{N} + \sum_{i=1}^{N-2} \frac{v_i}{N}) \\
    x(\frac{2v_{N-1}}{N} + \sum_{i=1}^{N-2} \frac{v_i}{N}) \\
    \vdots \\
    x(\frac{2v_{N-1}}{N} + \sum_{i=1}^{N-2} \frac{v_i}{N})
\end{bmatrix}.
\end{align*}
\end{proof}
We next study the quantity
$$C_N \colonequals \max_{2 \leq i \leq N - 1} \abs{\lambda_i}.$$ 
This constant $0 < C_N < 1$ describes the rate of convergence of a generic $N$-gon to its limit point. Indeed, working over $\C$, the rate of convergence of a polygon $v$ to the limit point under cyclic snap is invariant under the choice of basis on $\C^{N-1}$. For fixed $N$, working in the eigenbasis $w_{\lambda_1}, \ldots, w_{\lambda_{N-1}}$, we have
$$D^k w \ = \  \sum_{i = 1}^{N - 1} \lambda_i^k w_{\lambda_i}.$$ Since the simple eigenvalue $\lambda_1 = 1$ dominates, we have $\lim_{k \to \infty} D^k w = w_1$, and the distance to the limit point at step $k$ is the magnitude of $D^k w - w_1$. As a function of $k$, this magnitude is bounded by a constant multiple of $C_N^k$. It is easy to see that, choosing the appropriate eigenvector, a convergence rate of precisely $C_N^k$ can be achieved (allowing complex polygons).

Our last result shows that the rate of convergence asymptotically slows as $N$ increases.

\begin{prop} \label{prop_rate}
We have $\lim_{N \to \infty} C_N \ = \  1$.
\end{prop}

\begin{proof}
    Fix $N \geq 3$. If $P(\lambda) \ = \ 0$, then $Q(\lambda) \ = \  2\lambda^{N-1} - \lambda^{N-2} - 1 = 0.$
    Thus
    $$\lvert 2\lambda^{N-1} - \lambda^{N-2} \rvert \ = \  1,$$
    $$\lvert \lambda^{N-2} \cdot (2 \lambda-1 ) \rvert \ =\  1,$$
    $$\lvert \lambda \rvert^{N-2} \ =\  \frac{1}{\lvert 2\lambda - 1\rvert}.$$ 
Since $\lvert \lambda \rvert \le 1,$ we have $\lvert 2\lambda - 1 \rvert \le 3$ by the triangle inequality. Hence $\frac13 \le \lvert \lambda \rvert^{N-2}$, so $\left(\frac13\right)^{\frac{1}{N-2}} \leq \abs{\lambda} \leq 1$.
    This holds for all roots $\lambda$ of $P$, so taking the maximum,
    $$\left(\frac13\right)^{\frac{1}{N-2}} \ \le \  C_N \  \le \ 1.$$
    Then by the Squeeze Theorem,
    $$\lim\limits_{N \to \infty} \left( \frac13 \right)^{\frac{1}{N-2}} \ = \  \lim\limits_{N \to \infty} C_N \ = \  1.$$
\end{proof}
The same argument, taking a minimum over roots instead of a maximum, shows that the minimum absolute value $\abs{\lambda}$ also tends to $1$; in short, these root sets belong to annuli that tend to the unit circle. 

\section{Future Work}
This work examined the existence of a limiting point and the convergence rate $C_N$ to the limit point for one dynamical system, cyclic snap, defined by a rule $v'_i = \frac{1}{2}(v_{i-1} + v_{i+1})$ applied in a cyclic order to the vertices of a polygon. We collect here a few directions for future research.

\begin{enumerate}
    \item Rather than applying the snap rule in cyclic order, choose any finite sequence of indices taken from $0, \ldots, N - 1$ and apply the local snap operation in that composition order. This gives a range of dynamical systems; which ones admit a limit point, and is this limit point typically a rational linear combination of the initial vertices?
    \item Consider the randomized version of the above question: apply snap at a randomly chosen vertex at each time step. What is the probability of convergence to a limit point, and what is the expected rate of convergence?
    \item Are there any polygons that are periodic under cyclic snap, up to an affine transformation?
    \item What happens if one replaces the indices $i - 1$ and $i + 1$ in the snap rule with other shifts of $i$? For instance, what happens if we say instead that $v'_i = \frac{1}{2}(v_{i + a} + v_{i + b})$ for some choices of $a$ and $b$?
    \item Study the rate of convergence of the area to $0$ for $N$-gons.
    \item We suspect that the convergence rate $C_N$ increases in $N$, but we did not find a proof. Intuitively, it seems reasonable to guess that increasing the number of vertices causes more lag time for the convergence phenomenon to take effect.
\end{enumerate}

\section{Code}
To visualize polygon dynamics, we extended the code of Gil'i Zaid \cite{zaid_quadrilateral_folding}. The source code for our extended Streamlit applet is available in the \href{https://github.com/uihyeon123/Polymath-Jr.-Polygon-Dynamics-}{Polymath Jr. Polygon Dynamics repository} \cite{polymath_polygon_dynamics}. The entry point is \path{polygondynamics_2026-11-13.py}, which can be run locally using Streamlit. The applet supports polygons with four to seven vertices and includes diagonal reflection, reverse diagonal reflection, snap, recut, and fold followed by recut. It provides orbit plots, animations, and PNG and PDF figure downloads.

\begin{figure}[H]
\centering
\setlength{\fboxsep}{0pt}
\setlength{\fboxrule}{0.6pt}
\begin{subfigure}{0.4\linewidth}
\centering
\fbox{\includegraphics[width=0.99\linewidth]{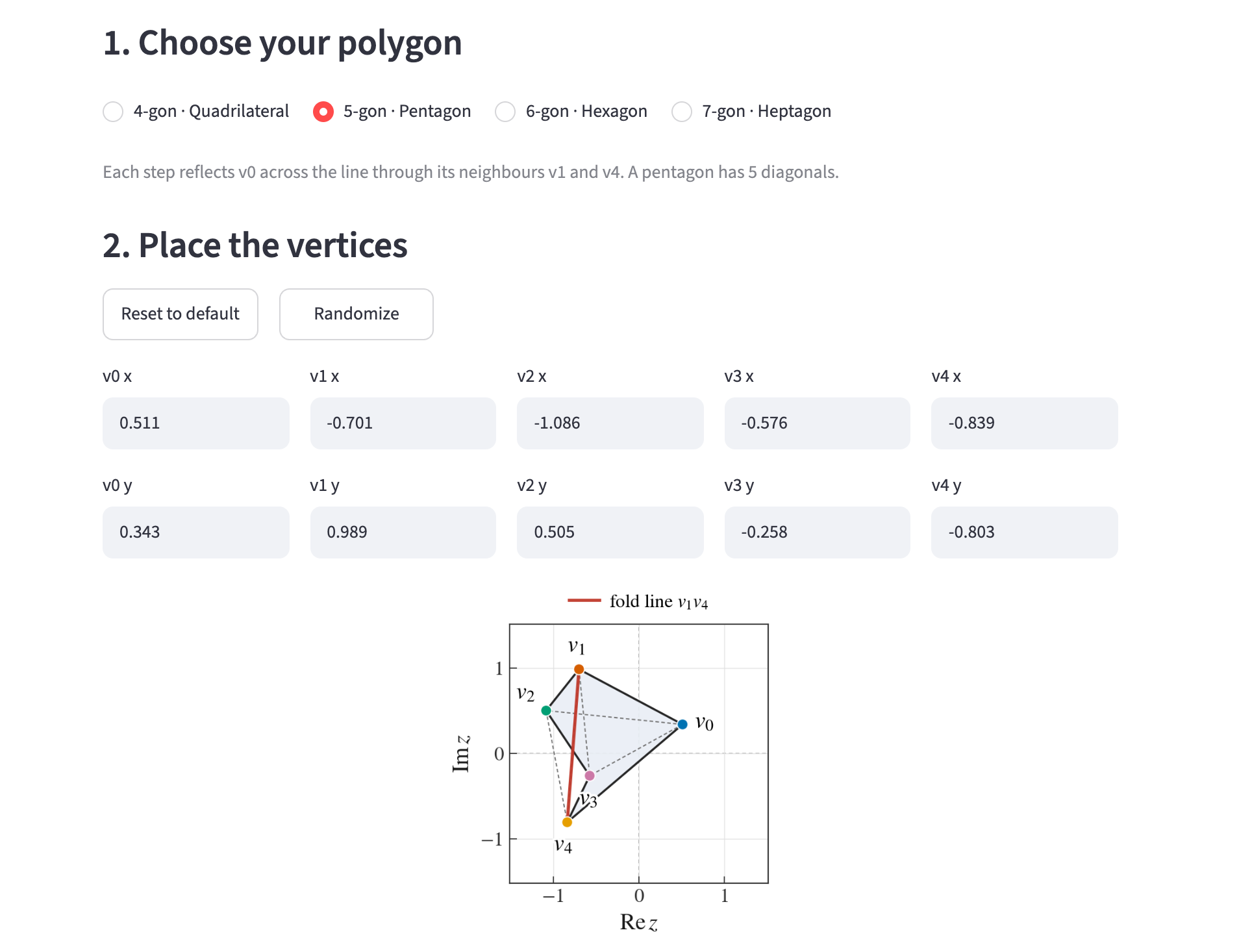}}
\caption{Polygon setup.}
\label{fig:code1}
\end{subfigure}

\medskip

\begin{subfigure}{0.5\linewidth}
\centering
\fbox{\includegraphics[width=0.99\linewidth]{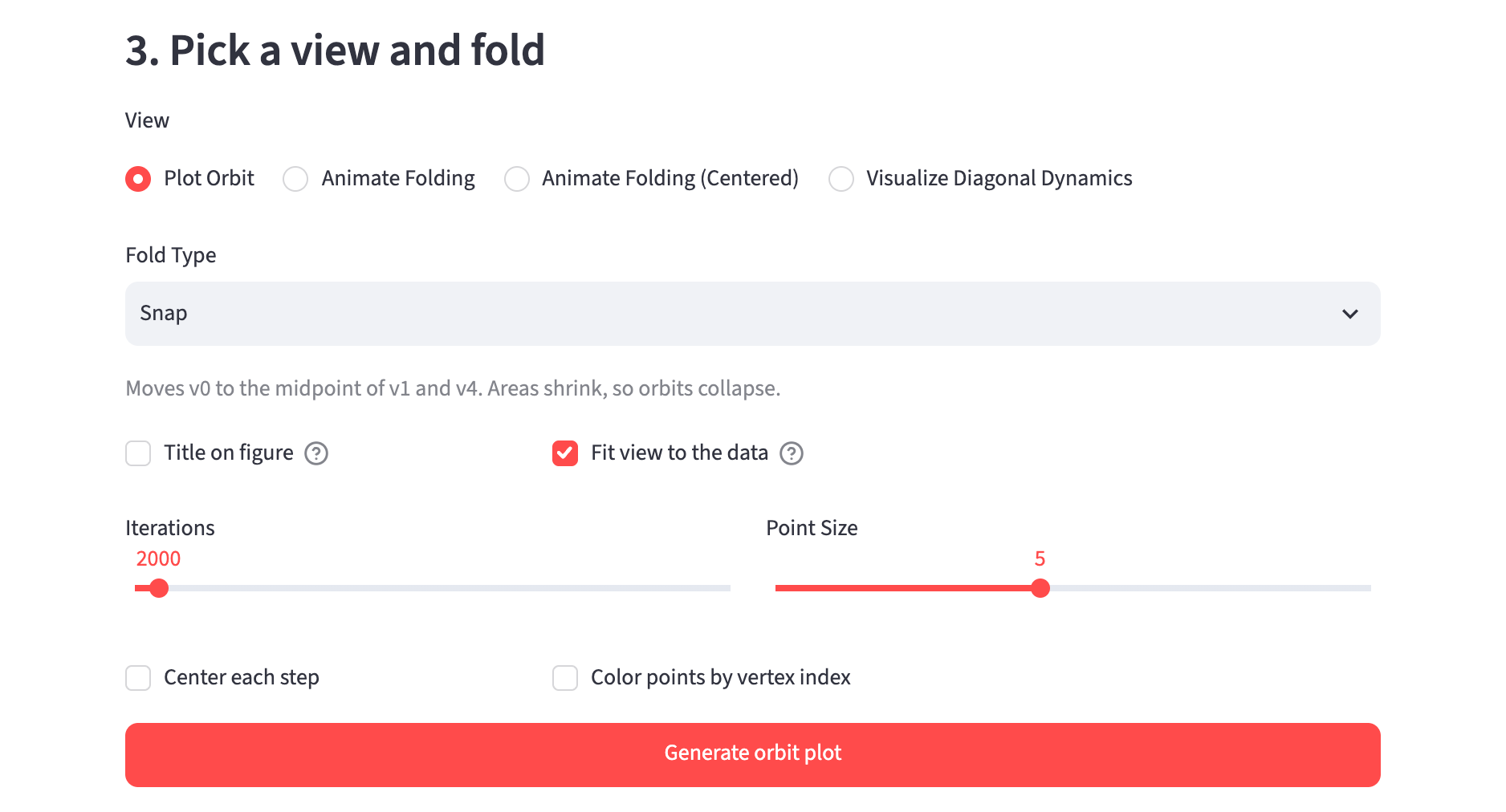}}
\caption{Dynamics settings.}
\label{fig:code2}
\end{subfigure}
\caption{The Streamlit interface.}
\label{fig:code}
\end{figure}

To explore snap, first choose the number of vertices and enter their coordinates under \textbf{Place the vertices}, as shown in Figure~\ref{fig:code1}. Select \textbf{Snap} under \textbf{Fold Type}, choose \textbf{Plot Orbit}, set the number of iterations, and click \textbf{Generate orbit plot}; see Figure~\ref{fig:code2}. Alternatively, choose \textbf{Animate Folding} and click \textbf{Generate animation} to view the successive polygons. Each applet iteration updates one vertex and cyclically relabels the vertices; thus $N$ applet iterations correspond to one application of cyclic snap $T$ on an $N$-gon.

\bibliographystyle{plain}

\bibliography{ref}
\end{document}